\documentclass[12pt,leqno]{article}

\usepackage{amssymb}
\usepackage{amsmath}
\usepackage{graphicx,psfrag,color}
\usepackage{url}
\usepackage{pstool}

\begin{document}

\title{Varieties generated by completions}
\author{Andr\'eka, H.\ and N\'emeti, I.}%
\date{}
\maketitle

\newcommand{\comment}[1]{}


\newcommand{\dm}[1]{\mathrm{dm(#1)}}
\newcommand{\rg}[1]{\mathrm{rg(#1)}}

\newcommand{\conv}{{}^{\smallsmile}}
\newcommand{\inv}{{}^{-1}}
\newcommand{\id}{\mbox{1'}}
\newcommand{\rp}{\mathrel{\vert}}
\newcommand{\Pow}{\mathcal{P}}
\newcommand{\E}{\mathcal{E}}

\newcommand{\di}{\mbox{\sf Di}}
\newcommand{\Z}{\mathcal{S}}
\newcommand{\HH}{\mathcal{H}}
\newcommand{\rep}{\mbox{\sf rep}}
\newcommand{\bb}{b^*}

\newcommand{\dow}{\mbox{$\downarrow\!\!(e,D)$}}

\newcommand{\de}{\mbox{\ :=\ }}
\newcommand{\deiff}{\mbox{$\quad : \Leftrightarrow\quad$}}
\newcommand{\Aa}{\mbox{$\mathfrak A$}}
\newcommand{\Bb}{\mbox{$\mathfrak B$}}
\newcommand{\Cc}{\mbox{$\mathfrak C$}}
\newcommand{\Dd}{\mbox{$\mathfrak D$}}
\newcommand{\Ee}{\mbox{$\mathfrak E$}}
\newcommand{\Mm}{\mbox{$\mathfrak M$}}
\newcommand{\qed}{\hfill\mbox{$\Box$}\bigskip}
\newcommand{\rr}{\mbox{\sf rep}}

\newcommand{\Ff}{\mbox{$\mathfrak F$}}
\newcommand{\Fff}{\mbox{$\Ff\times\Ff$}}
\newcommand{\Ffff}{\mbox{$\Ff\times \Ff\times \Ff$}}
\newcommand{\Gg}{\mbox{$\mathfrak G$}}
\newcommand{\Pp}{\mbox{$\mathfrak P$}}
\newcommand{\Rr}{\mbox{$\mathfrak R$}}

\newcommand{\G}{\mbox{$\mathcal{G}$}}

\newcommand{\RA}{\mbox{\sf RA}}
\newcommand{\SA}{\mbox{\sf SA}}
\newcommand{\RRA}{\mbox{\sf RRA}}
\newcommand{\CRA}{\mbox{\sf CRA}}
\newcommand{\V}{\mbox{\sf V}}
\newcommand{\W}{\mbox{\sf W}}
\newcommand{\PG}{\mbox{\sf PG}}
\newcommand{\LPG}{\mbox{\sf LPG}}
\newcommand{\RPG}{\mbox{\sf RPG}}
\newcommand{\K}{\mbox{\sf K}}
\newcommand{\Lo}{\mbox{\sf L}}

\newcommand{\Cm}{\mbox{\bf Cm}}
\newcommand{\Sub}{\mbox{\bf S}}
\newcommand{\Hom}{\mbox{\bf H}}
\newcommand{\Pro}{\mbox{\bf P}}
\newcommand{\Var}{\mbox{\sf Var}}

\newcommand{\scir}{;}

\newtheorem{thm}{Theorem}[section]
\newtheorem{lem}{Lemma}[section]
\newtheorem{Def}{Definition}[section]
\newtheorem{rem}{Remark}[section]
\newtheorem{prb}{Problem}[section]
\newtheorem{cor}{Corollary}[section]

\newcommand{\la}{\lambda}

\begin{abstract} We prove that persistently finite
algebras are not created by completions of algebras, in any ordered
discriminator variety. A persistently finite algebra is one without
infinite simple extensions. We prove that finite measurable relation
algebras are all persistently finite. An application of these
theorems is that the variety generated by the completions of
representable relation algebras does not contain all relation
algebras. This answers Problem 1.1(1) from  R. Maddux \cite{maddux}
in the negative. At the same time, we confirm the suggestion in that
paper that the finite maximal relation algebras constructed in M.
Frias and R. Maddux \cite{FM} are not in the variety generated by
the completions of representable relation algebras. We prove that
there are continuum many varieties between the variety generated by
the completions of representable relation algebras and the variety
of relation algebras.
\end{abstract}

\section{Introduction}\label{intro-sec}

Completions of partially ordered sets are obtained, intuitively, by
filling-in non-existent suprema. Various kinds of completions are in
use, for example in a join-completion we fill-in all suprema, and in
an ideal-completion we fill-in the suprema of nonempty directed sets
only. For Boolean algebras, most of these various kinds of
completions coincide with the so-called Dedekind-MacNeille
completion that we will simply call completion. Completions were
generalized from Boolean algebras to Boolean algebras with operators
by J. D. Monk \cite{monk}, where he also showed that the completion
of a relation algebra is again a relation algebra. It was natural to
expect that the completion of a representable relation algebra is
again representable. To a great surprise, this was disproved by I.
M. Hodkinson \cite{hod}. Maddux \cite{maddux} began to investigate
what kind of structures can a non-representable completion bring in.
He exhibits infinitely many finite non-representable relation
algebras that can be embedded into a completion of a representable
relation algebra. He then asks if the variety $\Var(\RRA^c)$
generated by all the completions of representable relation algebras
contains all relation algebras \cite[Problem1.1(1)]{maddux}.

In this paper, we deal with structures that cannot be created by
completions. We prove that in ordered discriminator varieties,
completions cannot create persistently finite algebras (Theorem
\ref{discr-thm}). In Frias-Maddux \cite{FM}, an infinity of maximal,
finite, integral relation algebras are constructed. Since these are
all persistently finite, they are all outside of $\Var(\RRA^c)$.
With this we confirm the suggestion in \cite{maddux} that these
might be outside of $\Var(\RRA^c)$. In this paper, we present
another sequence of persistently finite relation algebras that are
not integral. Namely, we prove that being finite measurable is a
persistent property (Theorem \ref{finmeas-thm}). All this gives a
negative answer to Problem 1.1(1) in \cite{maddux}. More can be
proved: there are continuum many varieties between $\Var(\RRA^c)$
and the variety $\RA$ of all relation algebras (Theorem
\ref{mad-thm}).

Section \ref{discr-sec} deals with discriminator varieties. We
define a general notion of density, and we prove that in ordered
discriminator varieties, dense extensions do not create new
persistently finite algebras. Sections \ref{pers-sec} and
\ref{compl-sec} deal with relation algebras. In section
\ref{pers-sec} we prove that finite measurable relation algebras are
all persistently finite, and we list some consequences of this.
In section \ref{compl-sec}, we use the theorems proved in the
previous sections to show that, in a sense, there is much room
between $\Var(\RRA^c)$ and $\RA$: there are continuum many varieties
between the two. However, this leaves open the question of how far
these varieties are from each other, that is whether $\Var(\RRA^c)$
is finitely axiomatizable over $\RA$ (this is \cite[Problem
1.1(2)]{maddux}).

\section{Dense extensions in discriminator
varieties}\label{discr-sec}

We recall the notion of a discriminator variety from \cite{ajn}. All
what we need in the present paper from this fascinating branch of
universal algebra is contained in \cite[section 1]{ajn} (a summary
which is taken from Werner \cite{werner}). A variety $\V$ is called
a \emph{discriminator variety} if it is generated by a class $\K$ of
algebras as a variety, and there is a term $\gamma$ in the language
of the variety that represents the \emph{quaternary discriminator}
function $g$ on each member of $\K$, where the function $g$ on a set
$A$ is defined by the property that
\[ g(x,y,u,v)=u \mbox{ when }x=y\quad\mbox{ and }\quad g(x,y,u,v)=v \mbox{ when }x\ne
y\] for all $x,y,u,v\in A$. We will not distinguish the function $g$
from the term $\gamma$ that defines it.

In this section we deal with ordered discriminator varieties as
defined in \cite[section 4]{ajn}. This is a discriminator variety
where a partial order $\le$ is defined by a finite set of equations:
\[ x\le y\qquad\mbox{ iff }\qquad \tau_i(x,y)=\sigma_i(x,y)\mbox{ for }i\le n \]
for some number $n$ and terms $\tau_i,\sigma_i,\ i\le n$ in the
language of the variety.
Discriminator varieties of Boolean algebras with operators are
ordered discriminator varieties where $\le$ is defined as the
Boolean ordering, in particular all varieties of relation algebras
or all varieties of finite dimensional cylindric algebras are
ordered discriminator varieties.

Since in ordered discriminator varieties the ordering can be an
arbitrary partial order, there is a great variety for the notion of
a completion (see, e.g., \cite{np97}). We will use a general notion
of density that covers most of completions.

\begin{Def}\label{dens-def}\label{dens-lem}
We say that $X\subseteq A$ is \emph{dense in $\langle A,\le\rangle$}
when for all $a\in A$ there is $x\in X$ with $x\le a$, and for all
$x\in X, x<a$ there is $y\in X$ with $x<y\le a$. \end{Def} This
density property is stronger than join-density (which means that
each element of $A$ is the supremum of the elements of $X$ below
it), but weaker than ideal-density (which means that each element of
$A$ is the supremum of its downset intersected with $X$ that is in
turn an ideal, i.e., a nonempty directed downward closed set). In
particular, a Boolean algebra is dense in the above sense in its
completion, see \cite[Lemma 15.1]{givbook2} or Lemma \ref{ddef-lem}.

Usually we denote the universes of algebras $\Aa, \Bb, \dots$ with
$A, B, \dots$.

We say that an ordered algebra $\Bb$ is a \emph{dense extension} of
$\Aa$, or that $\Aa$ is a \emph{dense subalgebra} of $\Bb$, when $A$
is dense in $\Bb$.
For a class $\K$ of similar algebras, $\K^d$ denotes the class of
dense extensions of members of $\K$ and $\Var\K$ denotes the variety
generated by $\K$.

Let $\Aa$ be an algebra and $\K$ a class of similar algebras. We say
that $\Aa$ is \emph{persistently finite in $\K$} if $\Aa$ is finite,
and all its simple extensions that are in $\K$ are finite.

The following theorem says that dense extensions do not generate new
persistently finite algebras, in an ordered discriminator variety.

\begin{thm}\label{discr-thm}
Assume that $\V$ is an ordered discriminator variety of finite
similarity type. Assume that $\Aa$ is simple and persistently finite
in $\V$. Then $\Aa\in\Var(\K^d\cap\V)$ implies that $\Aa\in\Var\K$,
for any subclass $\K$ of $\V$.
\end{thm}

\noindent{\bf Proof.} Let $\V$, $\K$ and $\Aa$ be as in the
statement of the theorem. Assume that $\Aa\in\Var(\K^d\cap \V)$, we
will show that $\Aa\in\Var\K$.
Let $\Hom\Lo$, $\Sub\Lo$ and $\Pro\Lo$ denote the classes of all
homomorphic images, all subalgebras and all direct products,
respectively, of members of $\Lo$, for $\Lo$ a class of similar
algebras. Thus, $\Var\Lo=\Hom\Sub\Pro\,\Lo$, for any class $\Lo$.
If $\Dd_i$ is a dense subalgebra of $\Ee_i$ for all $i\in I$ then
the direct product of the $\Dd_i, i\in I$ is a dense subalgebra of
the direct product of the $\Ee_i, i\in I$, this is straightforward
to check. Thus we have
\[\tag{1}\label{pro} \Pro(\K^d\cap\V)\subseteq(\Pro\K)^d\cap\V .\]
By $\Aa\in\Var(\K^d\cap\V)$ and \eqref{pro} we have that
$\Aa\in\Hom\Sub((\Pro\K)^d\cap\V)$, so there are algebras $\Bb, \Ee,
\Dd$ and a homomorphism $h$ such that
\begin{description}
\item{$\bullet$} $h:\Bb\to\Aa$ is a surjective homomorphism,
\item{$\bullet$} $\Bb$ is a subalgebra of $\Ee\in\V$,
\item{$\bullet$} $\Dd$ is a dense subalgebra of $\Ee$\quad and\quad $\Dd\in\Pro\K$.
\end{description}
Note that $\Ee\notin\Var\K$ may be the case. We will show that
$\Aa\in\Sub\Hom\{\Dd\}\subseteq\Var\K$.

Let $a_1,\dots, a_n$ be a repetition-free listing of the elements of
$A$. There is such a listing because $\Aa$ is finite. Since $\Aa$ is
simple, it has at least two elements, so $n\ge 2$. Let
$b_1,\dots,b_n\in B$ be such that $h(b_j)=a_j$ for all $1\le j\le
n$. There are such elements because $h$ is surjective. We may assume
that $\Bb$ is generated by $\{ b_1,\dots,b_n\}$, because of the
following. $\Aa$ is a homomorphic image of the subalgebra $\Bb'$ of
$\Bb$ generated by $b_1,\dots,b_n$ and $\Bb'$ is a subalgebra of
$\Ee$ by $\Bb'\subseteq\Bb\subseteq\Ee$, so we can choose this
$\Bb'$ in the first step. In the following, assume that $\Bb$ is
generated by $b_1,\dots,b_n$.

By \cite[Lemma 2.1]{ajn}, the kernel $\ker(h)$ of $h$ is a compact
congruence of $\Bb$. Since $\V$ is a discriminator variety,
compactness implies that $\ker(h)$ is principal, see
\cite[Thm.1.1.III.(iii)]{ajn}. Let 
$r,s\in B$ be such that
\[\tag{3}\label{3}\mbox{$\ker(h)$ is generated by the pair
$\{\langle r,s\rangle\}$ as a congruence in $\Bb$.}\]

\noindent We want to use the property of $\Aa$ that all its simple
extensions in $\V$ are finite. To this end, let us take a subdirect
decomposition of $\Ee$ to simple factors, say, $\Ee$ is a subdirect
product of $\langle\Cc_i : i\in I\rangle$ with projections
$\pi_i:\Ee\to\Cc_i$ where each $\Cc_i$ is simple. There is such a
decomposition by Birkhoff's subdirect decomposition theorem and
because in discriminator varieties the subdirectly irreducible
algebras are exactly the simple ones, see \cite[Thm.1.1.II]{ajn}.
\noindent Let
\[\mbox{$I_0=\{ i\in I : \pi_i(r)=\pi_i(s)\mbox{ and }\pi_i(b_1)\ne\pi_i(b_2) \}$.}\]
We show that
\[\tag{4}\label{4} \mbox{$h_i=\{\langle h(b),\pi_i(b)\rangle : b\in B\}$ embeds $\Aa$ into
$\pi_i[\Bb]\subseteq\Cc_i$, for all $i\in I_0$.}\]
Indeed, let $i\in I_0$. Then $\langle r,s\rangle\in\ker(\pi_i)$ by
the definition of $I_0$. Thus $\ker(h)\subseteq\ker(\pi_i)$ because
$\ker(h)$ is generated by this pair as a congruence in
$\Bb\subseteq\Ee$. Now, $\pi_i:\Bb\to\Cc_i$ since $\Bb\subseteq\Ee$
and $\pi_i:\Ee\to\Cc_i$. Therefore $h_i:\Aa\to\pi_i[\Bb]$, by
$\ker(h)\subseteq\ker\pi_i$ and $h:\Bb\to\Aa$. Since $\Aa$ is
simple, then either $h_i$ is an embedding, or $h_i$ maps $\Aa$ onto
the one-element algebra. However, this latter case would entail
$\pi_i(b_1)=\pi_i(b_2)$, which is not the case by the definition of
$I_0$ and $i\in I_0$. Thus $h_i$ is an embedding, and the proof of
\eqref{4} is complete.

Now, by using persistently finiteness of $\Aa$, we show that there
is a natural number $k$ such that
\[\tag{5}\label{5}\mbox{the cardinality of $\Cc_i$ is smaller than $k$, for all $i\in I_0$.}\]
Indeed, $\Cc_i$ is simple and $\Cc_i\in\V$ for all $i\in I$ by our
assumption $\Ee\in\V$. Then $h_i:\Aa\to\Cc_i$ implies that $\Cc_i$
is finite. However, we have to show more, we have to show the
existence of the finite upper bound $k$. Proving by contradiction,
assume that $\Aa$ has ever larger simple extensions $\Mm_n$ in $\V$.
Then $\Aa$ is embeddable into an infinite ultraproduct  $\Mm$ of
these. Now, $\Mm\in\V$ because $\V$ is a variety, so it is closed
under taking ultraproducts. Since $\V$ is a discriminator variety,
the class of its simple members is closed under taking ultraproducts
(\cite[Thm.1.1.II.(iv)]{ajn}), and hence $\Mm$ is simple. This
contradicts persistently finiteness of $\Aa$, and so there is a
finite upper bound $k$ for the simple $\V$-extensions of $\Aa$. This
proves (5).

Next we show that \[\tag{6}\label{5a}\mbox{$I_0\ne\emptyset$.}\]
Indeed, let $\bb=g(r,s,b_1,b_2)$, where $g$ is the quaternary
discriminator term. Now, $\bb\in B$ by its definition, and
$h(\bb)=
g(h(r),h(s),h(b_1),h(b_2))$ because $h$ is a homomorphism  and $g$
is a term in the language. Since $\ker(h)$ is generated by $\langle
r,s\rangle$, we have that $h(r)=h(s)$. Thus, $h(\bb)=h(b_1)=a_1$, by
the definition of $g$. Assume that $i\notin I_0$. Then either
$\pi_i(r)\ne\pi_i(s)$ or $\pi_i(b_1)=\pi_i(b_2)$. In either case,
$\pi_i(\bb)=\pi_i(b_2)$. Assume that $I_0=\emptyset$, then
$\pi_i(\bb)=\pi_i(b_2)$ for all $i\in I$, which implies $\bb=b_2$.
This is a contradiction, since $h(\bb)=a_1$ while $h(b_2)=a_2\ne
a_1$. Thus, $I_0\ne\emptyset$ and (6) has been proved.

We begin a ``transit" from $\Ee$ to $\Dd$. So far, (4) implies that
$h_i[\Aa]\subseteq\pi_i[\Bb]\subseteq\pi_i[\Ee]$. Note that $\Dd$ is
a subalgebra of $\Ee$ and $\Dd\in\Pro\K$. Our plan is to find an
$i\in I_0$ such that $h_i[\Aa]\subseteq\pi_i[\Dd]$.  This would show
that $\Aa\in\Sub\Hom\{\Dd\}\subseteq\Var\K$.
To this end, we want to find $d_1,\dots,d_n\in D$ such that there is
an $i\in I_0$ such that $\pi_i(b_1)=\pi_i(d_1),\dots,
\pi_i(b_n)=\pi_i(d_n)$. This would imply that
$\pi_i[\Bb]\subseteq\pi_i[\Dd]$ because $\Bb$ is generated by
$b_1,\dots,b_n$, and then \eqref{4} implies that
$h_i[\Aa]\subseteq\pi_i[\Dd]$.

From now on, the following notation will be convenient to use. Let
$\alpha(e_1,\dots,e_m)$ be a conjunction of equations and
non-equations of terms of elements from $E$. These are all open Horn
formulas.
$I(\alpha)$ denotes those factors from $I$ where $\alpha$ is true
under taking the projections:
\[ I(\alpha(e_1,\dots,e_m)) = \{ i\in I :
\Cc_i\models\alpha(\pi_i(e_1),\dots,\pi_i(e_m))\} . \] With this
notation, $I_0=I(r=s \land b_1\ne b_2)$.
We say that $J\subseteq I$ is \emph{definable} if $J$ is
$I(\alpha(e_1,\dots,e_m))$ for some open Horn formula $\alpha$ and
elements $e_1,\dots,e_m$ of $E$. In discriminator varieties, to each
open Horn formula $\alpha$ there is an equation $e$ such that in
each simple member of the variety, $\alpha$ and $e$ have the same
truth-evaluations, see \cite[Thm.1.1.V.]{ajn}.
Therefore, there are terms $\rho,\delta$  such that
$I(\alpha(e_1,\dots,e_m)) =
I(\rho(e_1,\dots,e_m)=\delta(e_1,\dots,e_m))$. Since
$\rho(e_1,\dots,e_m), \delta(e_1,\dots,e_m)$ are elements of $E$, we
have that $J$ is definable exactly when there are $p,q\in E$ such
that $J=\{ i\in I : \pi_i(p)=\pi_i(q)\}$. In this case, we say that
$J$ is defined by $p=q$. Assume that $J\subseteq I$. By
``$\alpha(e_1,\dots,e_m)$ holds on $J$" we understand that
$\alpha(\pi_i(e_1),\dots,\pi_i(e_m))$ holds for all $i\in J$. We
note that in ordered discriminator varieties, $\tau<\sigma$ counts
as an equation, because $<$ is defined by a finite set of equations
which is an open Horn formula.

We say that $J\subseteq I$ is \emph{bounded} if there is a natural
number $k$ such that $|C_i|\le k$ for all $i\in J$. Finally, when
$a,b\in E$ and $J\subseteq I$, we say that $a$ \emph{approximates}
$b$ in $J$ if there is $j\in J$ such that $a_j=b_j$. To find our
$d_1,\dots,d_n$, we will use the following statement (7) repeatedly.
We note that this statement is the ``heart" of the proof of Theorem
\ref{discr-thm}.

\begin{description}
\item{(7)} Each element of $E$ can be approximated by an element of $D$ in any
nonempty definable bounded subset of $I$.
\end{description}
Indeed, to prove (7), let $b\in E$ and let $J$ be a nonempty
definable bounded subset of $I$. By density of $D$ in $E$, there is
an $a\in D$ such that $a\le b$. If $a$ approximates $b$ in $J$ then
we are done. Assume that $a$ does not approximate $b$ in $J$, then
\[\tag{8} \mbox{$a<b$\ \  on\ \  $J$.}\]
Define $b_1=g(p,q,b,a)$, where $J$ is defined by $p=q$. Then $b_1$
is $b$ on $J$ and $b_1$ is $a$ outside $J$ and $a<b_1$ by (8).
Choose $a_1\in D$ such that $a<a_1\le b_1$, there is such an $a_1$
by Definition \ref{dens-def}. If $a_1$ approximates $b$ on $J$, then
we are done. So, assume that $a_1$ does not approximate $b$ on $J$
and let $J_1=I(a<a_1)$. Then $J_1$ is a subset of $J$ because
$a<a_1\le b_1$ but $a$ and $b_1$ agree outside of $J$. (We note that
the purpose of using $b_1$ in place of $b$ was to achieve
$J_1\subseteq J$.) By $a<a_1$ we have that $J_1$ is nonempty. $J_1$
still has the bound $k$ because it is a subset of $J$. Also, $J_1$
is definable by $J_1=I(a<a_1)$. Thus, $J_1$ is a nonempty, definable
bounded subset of $J$ and
\[ \mbox{$a<a_1<b$\ \  on\ \  $J_1$.}\]
We proceed this way: assume that for $m$ we already defined
$a_1,\dots,a_m\in D$ and a nonempty definable subset $J_m$ of $J$
such that
\[ \mbox{$a<a_1<\dots<a_m<b$\ \ on\ \ $J_m$.}\]
Having this, by $a_m<b$ and Definition \ref{dens-def} we can find an
$a_{m+1}\in D$ such that $a_m<a_{m+1}<b$. If $a_{m+1}$ approximates
$b$ on $J_m$ then we are done. If not, then let
$b_{m+1}=g(p_m,q_m,b,a_{m+1})$ where $J_m$ is defined by $p_m=q_m$.
Let $J_{m+1}=I(a_m<a_{m+1})$, this is a nonempty definable subset of
$J_m$ such that
\[\tag{9} \mbox{$a<a_1<\dots<a_m<a_{m+1}<b$\ \ on\ \ $J_{m+1}$.}\]
In particular, (9) implies that $|C_i|\ge m+3$ for all $i\in
J_{m+1}$. Because $J$ is bounded, this process cannot be continued
ad infinitum, so at one of the steps we have to find an
approximation of $b$ as desired. This proves (7).

We turn to finding $d_1,\dots,d_n\in D$ such that
$I(b_1=d_1\land\dots\land b_n=d_n)\cap I_0\ne\emptyset$. We have
seen that $I_0$ is a nonempty definable bounded subset of $I$, see
(5), (6). Assume that $p=q$ defines $I_0$. We begin with $b_1$. By
(7), there is $d_1\in D$ which approximates $b_1$ on $I_0$. Let
$I_1=I(b_1=d_1\land p=q)$. Then $I_1$ is a nonempty definable
bounded subset of $I_0$. By (7), there is $d_2\in D$ which
approximates $b_2$ in $I_1$. Let $I_2=I(b_1=d_1\land b_2=d_2\land
p=q)$. Then $I_2$ is a nonempty definable bounded subset of $I_0$,
and so on. In the last step we get a $d_n\in D$ such that $d_n$
approximates $b_n$ on $I_{n-1}=I(b_1=d_1\land\dots\land
b_{n-1}=d_{n-1}\land p=q)$. Let $J=I(b_1=d_1\land\dots\land
b_{n}=d_{n}\land p=q)$. Then $J$ is nonempty, let $i\in J$ be
arbitrary. Then, $\pi_i(b_1),\dots,\pi_i(b_n)\in\pi_i[\Dd]$ by
$d_1,\dots,d_n\in D$ and
$\pi_i(b_1)=\pi_i(d_1)\dots,\pi_i(b_n)=\pi_i(d_n)$. By (4) we have
that $h_i(a_j)=\pi_i(b_j)$ for all $1\le j\le n$, thus
$h_i[\Aa]\subseteq\pi_i[\Dd]$ because $A=\{ a_1,\dots,a_n\}$. By
assumption we have that $\Dd\in\Pro\K$, so
$\Aa\in\Sub\Hom\Pro\,\K\subseteq\Var\K$. With this, the proof of
Theorem \ref{discr-thm} is complete. \qed

\section{Persistently finite non-representable relation
algebras}\label{pers-sec}

An algebra $\Aa=\langle A,+,-,;,\conv,\id\rangle$ is a
\emph{concrete algebra of binary relations} if $\Aa$ is a set of
binary relations with a biggest one, and the operations
$+,-,;,\conv,\id$ are, respectively, the following natural
operations on binary relations: union of two relations, taking the
complement of a relation with respect to the biggest relation,
relational composition of two relations, converse of a relation, and
the identity relation on the domain of the biggest relation. The
class of all algebras isomorphic to concrete algebras of binary
relations is denoted by $\RRA$. This is a variety which is not
definable by a finite set of equations (classic results due to A.
Tarski and J. D. Monk, respectively). The variety $\RA\supseteq\RRA$
of relation algebras is a finitely axiomatized variety that
approximates $\RRA$ surprisingly well: an algebra $\Aa=\langle
A,+,-,;,\conv,\id\rangle$ is a \emph{relation algebra} if $\langle
A,;,\conv,\id\rangle$ is an involuted monoid,  $\langle
A,+,-,;,\conv\rangle$ is a Boolean algebra with normal and additive
operators,  and one more identity true of concrete algebras of
binary relations also holds in it, namely $r\conv;-(r;s)\le-s$. We
use $\le, 0,1,\cdot$ with their usual definitions in a Boolean
algebra. The elements of $\RRA$ are called representable relation
algebras and the elements of $\RA\setminus\RRA$ are called
\emph{non-representable} relation algebras. If $\Aa\in\RA$, we
assume that its operations are as above.

Assume that $\Aa\in\RA$ and $y\in A$. We say that $y$ is a
\emph{functional element} if $y\conv;y\le\id$.  A relation algebra
$\Aa$ is called \emph{measurable}, if the identity constant $\id$ is
the supremum of atoms, and each atom $x\le\id$ is \emph{measurable}
in the sense that $x;1;x$ is the supremum of the functional elements
below it. The number of the functional elements below $x;1;x$ is
called the \emph{measure of $x$}. These names reflect their meanings
in concrete algebras of binary relations. Namely, assume $\Aa$ is
such. Then $y\in A$ is functional exactly when $y$ is a function as
a relation, a subidentity element $x$ corresponds to a subset $X$ of
the domain of the biggest relation via $x=\{\langle u,u\rangle :
u\in X\}$. When the biggest element of $\Aa$ is of form $U\times U$
then $x;1;x$ is just the square $X\times X$ and it can be showed
that the measure of a measurable atom $x$ coincides with the size
$|X|$ of $X$.

The following theorem says that the property of being finite and
measurable is persistent in $\RA$, i.e., this property is preserved
by simple relation algebra extensions.

\begin{thm}\label{finmeas-thm}
Assume that $\Mm$ is a finite measurable relation algebra. If $\Mm$
can be embedded into a simple relation algebra $\Aa$, then $\Aa$
itself is finite and measurable.
\end{thm}

\noindent{\bf Proof.} Assume that $\Mm$ is a finite measurable
relation algebra and $\Mm\subseteq\Aa\in\RA$ where $\Aa$ is simple.
Let $I$ denote the set of the subidentity atoms--- atoms below the
identity $\id$--- of $\Mm$, and for all $x\in I$, let $F_x$ denote
the set of functional elements below $x\scir 1\scir x$. Then $I$ is
finite, and for all $x\in I$, the set $F_x$ is finite, too, because
$\Mm$ is finite. Further, $x;1;x$ is the sum of $F_x$. Temporarily,
let us fix $x\in I$.

By $\Mm\subseteq\Aa$, we have that $e=x\scir 1\scir x$ is an element
of $\Aa$, too, and it is an equivalence element, that is to say,
$e\conv=e=e\scir e$ because $\Aa\in\RA$, see \cite[Lemma
5.64]{givbook1}. Then the relativization $\Aa(e)$ of $\Aa$ to $e$ is
also a relation algebra \cite[Theorem 10.1]{givbook1}. It is simple
because $\Aa$ is simple and $x;1;x$ is a nonzero square
\cite[Theorem 10.8]{givbook1}.
We note that the universe of $\Aa(e)$ is the downset of $e$ in
$\Aa$, and the operations $+,;,\conv$ of $\Aa(e)$ are the same as
those of $\Aa$ while $-,\id$ in $\Aa(e)$ are the relativized
versions of those in $\Aa$ which means that $-(a)=e-a$ and the
identity of $\Aa(e)$ is $\id\cdot e$.

Each $g\in F_x$ is functional in $\Mm$, so $g$ is functional in
$\Aa$, and functional in $\Aa(e)$, too. In $\Mm$ we have that
$e=\sum F_x$, so this same is true in $\Aa$ and in $\Aa(e)$,  since
$F_x$ is finite. Thus, in $\Aa(e)$ the unit, $e$ is a sum of
finitely many functional elements.
Therefore, $\Aa(e)$ is representable on a finite set because it is
simple, by \cite[Theorem 4.32]{jt52}. This means that $\Aa(e)$ can
be represented such that the unit has the form of $U\times U$ for
some finite set $U$.
A finite relation has only finitely many subsets, so, $\Aa(e)$ is
finite and so atomic. Since subsets of functions are functions and
the unit in $\Aa(e)$ is the sum of finitely many functional
elements, we have that each element in $\Aa(e)$ is the sum of
functional elements. In particular, $\Aa(e)$ is measurable with
finitely many subidentity atoms of finite measure. Since the
universe of $\Aa(e)$ is the set of all elements of $A$ that are
below $e$, we get that in $\Aa$, too, the identity element $x$ of
$\Aa(e)$ is the sum of finitely many subidentity atoms of finite
measure.

A relation algebra is called \emph{finitely measurable} if it is
measurable and each subidentity atom has finite measure in it.
The identity element of $\Aa$ is the sum of $I$, the set of
subidentity atoms of $\Mm$, since this is true of $\Mm\subseteq\Aa$.
We have seen that each $x\in I$ in $\Aa$ is the finite sum of
measurable atoms of finite measure. So in $\Aa$, too, the identity
element $\id$ is the sum of finitely many atoms of finite measure.
Thus, $\Aa$ is finitely measurable. Then $\Aa$ is atomic by
\cite[Theorem 8.3]{gaapal}. We do not know yet whether it is
complete or not because we do not know yet that it is finite, since
in principle, there might be infinitely many atoms below $x;1;y$ for
distinct subidentity atoms $x,y$.

We now use the representation theorem Theorem 7.4(ii) from
\cite{gaapal}. It says that the completion $\Cc$ of $\Aa$ is a coset
relation algebra that is determined by a group coset frame
$(\G,\varphi,S)$ distilled from $\Aa$. We are going to show that
$\Cc$ is finite. The system $\G$ of groups is $\langle \Gg_x : x\in
J\rangle$ where $J$ is the set of subidentity atoms of $\Aa$ and
$\Gg_x$ is the group of functional elements below $x;1;x$ for $x\in
J$. We have already shown that $J$ is finite and $\Gg_x$ is finite
for all $x\in J$. Now, the atoms below $x;1;y$ in $\Cc$ are in
one-to-one correspondence with $\Gg_x\slash H_{xy}$ where $H_{xy}$
is a normal subgroup of $\Gg_x$, by the definition of an algebra
determined by $\G$ (see \cite[p.1171]{gaapal}).
Since $\Gg_x$ is finite, this means that there are finitely many
atoms below $x;1;y$, for any $x,y\in J$. Since $J$ is finite, this
means that $\Cc$ is finite. Since $\Aa$ is a subalgebra of $\Cc$, we
have that $\Aa$ is finite and we have already shown that $\Aa$ is
measurable. \qed

There are many finite and infinite representable measurable relation
algebras, these are described in \cite{giv1}. Also, infinitely many
finite and infinite non-representable measurable relation algebras
are constructed in \cite[sections 3, 4]{AGN18}. (We note that these
are not weakly representable, either.)
%
By using Theorem \ref{finmeas-thm}, the non-representable measurable
algebras can be used to give answers to Problems 2,3,4 from
\cite{ajn} (see also Problems P5,P6,P7 in \cite{madpers}). These
problems were already solved in \cite{FM}, but measurable algebras
provide a different kind of examples for their solutions. Below, we
elaborate on this.

A relation algebra is called \emph{maximal} if it is finite, simple
and has no proper simple relation algebra extension. It is known
that the representable maximal relation algebras are exactly the
finite full set relation algebras, that is to say, the concrete
algebras of all subsets of $U\times U$ for some finite set $U$.
Until 1997, these were the only known maximal relation algebras.

Problem 3 in \cite{ajn} asks whether there are any non-representable
simple absolute retracts in $\RA$ (or in the variety $\SA$ of
semi-associative relation algebras). In semi-simple varieties $\V$,
being an absolute retract and being maximal among the simple
algebras are equivalent \cite[Lemma 3.3]{ajn}. This problem was
solved in \cite{FM}, Frias and Maddux constructed infinitely many
non-representable maximal relation algebras. These algebras are
integral, i.e., the identity constant $\id$ is an atom in them.
Theorem \ref{finmeas-thm}, together with \cite[Theorem 4.2]{AGN18},
provide different, non-integral maximal examples for Problem 3, as
follows.

We have seen when proving (5) in the proof of Theorem
\ref{discr-thm} that a persistently finite algebra can have only
finitely many simple extensions (up to isomorphism, in a
discriminator variety of finite similarity type). Therefore, each
simple persistently finite algebra can be extended to a maximal one.
The maximal extension of a representable measurable relation algebra
is a full set relation algebra, this provides many persistently
finite and non-maximal simple relation algebras. On the other hand,
the maximal extension of a non-representable measurable relation
algebra is also non-representable, and thus provides new solutions
for \cite[Problem 3]{ajn}. The simple non-representable measurable
relation algebras are all non-integral, because in an integral
simple measurable relation algebra the square $\id;1;\id=1$, and
thus measurability implies that the unit is the supremum of
functional elements and these are known to be all representable.

Problem 4 in \cite{ajn} asks whether for an atom in a relation
algebra satisfying the equality $p;1;p\conv\le\id$ is necessary or
not for being persistent. We note that any maximal non-representable
relation algebra gives a negative answer for this problem, because a
finite algebra is atomic, if it is maximal, then all its atoms are
persistent, and it is known that if all atoms $p$ in an algebra
satisfy the given equality $p;1;p\conv\le\id$, then the algebra must
be representable.
This problem was solved in the negative in \cite{FM}, too, and the
maximal non-representable measurable relation algebras provide
further examples disproving an affirmative answer.

Problem 2 in \cite{ajn} asks whether there exists a simple relation
algebra that is not embeddable into a one-generated relation
algebra. Any maximal relation algebra that is not one-generated
gives a negative answer to this problem. Indeed, \cite{FM} gave such
integral examples. We believe that the smallest non-representable
measurable relation algebra constructed from the two-element group
in \cite[section 5]{andgiv1} and in \cite[sections 3,4]{AGN18} is
maximal, 2-generated but not one-generated. (However, we did not
check the details of this claim.) That would provide a new example
for a negative answer to \cite[Problem 2]{ajn}.

We note that Problem 5 from \cite{ajn} has been solved recently by
Mohamed Khaled. He proved that the finitely generated free
non-associative relation algebras are not atomic \cite{khal1}, and
the finitely generated finite dimensional free non-commutative
cylindric algebras are not atomic, either \cite{khal2}.
Problem 1 of \cite{ajn} asking whether the free $m$-generated and
free $n$-generated pairing algebras may be isomorphic for distinct
$m,n$ is still open, to our best knowledge. With this, we surveyed
the present-day statuses of all the problems given in \cite[section
11]{ajn}.

\section{Completions of representable relation algebras}\label{compl-sec}

In this section we are going to apply Theorem \ref{discr-thm} in the
context of relation algebras. We begin with showing that the notion
of a dense extension as defined below Definition \ref{dens-def} in
this paper is equivalent with the one most widely used in the
literature, for Boolean algebras with additional operators, and in
particular, for relation algebras. (See, for example,
\cite[Definition 15.4(ii)]{givbook2}, \cite[Definition 2.6]{hh02},
\cite[p.237]{ma06}.)
Let $\le$ denote the Boolean ordering, i.e., $x\le y$ is defined by
$x+y=y$.

\begin{lem}\label{ddef-lem}
Assume that $\Aa$ and $\Bb$ are Boolean algebras with operators and
$\Aa\subseteq\Bb$. Then {\rm (i)} and {\rm (ii)} below are
equivalent.
\begin{description}
\item[(i)] $\Aa$ is a dense subalgebra of $\Bb$.
\item[(ii)] For all nonzero $b\in B$ there is a nonzero $a\in A$
such that $a\le b$.
\end{description}
\end{lem}

\noindent{\bf Proof.} Assume (i) and let $0<b\in B$. Since
$\Aa\subseteq\Bb$, we have that $0\in A$. Then, $0\in A$, $0<b\in B$
and (i) imply, by Definition \ref{dens-def}, that there is $a\in A$
such that $0<a\le b$, and we are done.

Assume now (ii), we want to show that $A$ is dense in $\Bb$. Each
element of $B$ has an element of $A$ below it, namely $0\in A$.
Assume that $a<b\in B$ and $a\in A$. Then $0<b-a\in B$, so by (ii)
there is $x\in A$ such that $0<x\le b-a$. Then $a+x\in A$ and
$a<a+x\le b$, and we are done. \qed

A completion $\Rr^c$ of a relation algebra $\Rr$ is defined as a
complete, dense extension of $\Rr$; this exists and is unique up to
an isomorphism that leaves $R$ fixed. (See, for example,
\cite[Definition 15.17]{givbook2}, \cite[Definition 2.25, Lemma
2.26]{hh02}, \cite[p.323]{ma06}.)
For a class $\K$ of relation algebras, let $\K^c$ denote the class
of completions of elements of $\K$. The following lemma implies that
$\Var(\K^d)=\Var(\K^c)$ for any $\K\subseteq\RA$.

\begin{lem}\label{code-lem}
$\Sub\K^d=\Sub\K^c\subseteq\RA$ for any $\K\subseteq\RA$.
\end{lem}

\noindent{\bf Proof.} To show $\K^d\subseteq\Sub\K^c$, let
$\Aa\in\K$ and let $\Bb$ be a dense extension of $\Aa$. We want to
show that $\Bb\in\Sub\K^c$. Let $\Bb^c$ be a completion of $\Bb$,
this exists by $\K\subseteq\RA$. Then $\Bb^c$ is a dense extension
of $\Aa$ and it is complete, thus $\Bb^c$ is a completion of $\Aa$,
and therefore $\Bb\subseteq\Bb^c\in\K^c$.  The other direction
follows from the definition of a completion, this definition
immediately implies that $\K^c\subseteq\K^d$. Since $\RA$ is closed
under completions, we immediately get $\K^c\subseteq\RA^c=\RA$. \qed

It is well known that the variety $\RA$ of relation algebras is a
discriminator variety (see, for example, \cite[Corollary 5.7]{ajn},
or \cite[Theorem 3.19]{hh02}, \cite[p.386]{ma06}, \cite{jip}).
Since $\RA$ has a Boolean reduct, it is an ordered discriminator
variety with letting $\le$ be the Boolean order. Thus, Theorem
\ref{discr-thm} can be applied with taking $\V$ to be the variety of
relation algebras.
We are ready to state the theorem that gives an answer to
\cite[problem 1.1(1)]{maddux}.

\begin{thm}\label{mad-thm} $\Var(\RRA^c)\ne\RA$. Moreover, the following {\rm (i)-(iv)} hold.
\begin{description}
\item[(i)]
There are infinitely many finite simple integral relation algebras
which are not in $\Var(\RRA^c)$.
\item[(ii)]
There are infinitely many finite simple non-integral relation
algebras which are not in $\Var(\RRA^c)$.
\item[(iii)]
There are continuum many varieties $\W$ such that
$\Var(\RRA^c)\subset\W\subset\RA$.
\item[(iv)]
A finite measurable relation algebra is in $\Var(\RRA^c)$ exactly
when it is representable.
\end{description}
\end{thm}

\noindent{\bf Proof.}
First we prove (iv). Assume that $\Mm$ is a finite measurable
relation algebra. If $\Mm$ is representable, then it is its own
completion since it is complete, and so it is in $\RRA^c$. In the
reverse direction, assume that $\Mm\in\Var(\RRA^c)$, we will show
that $\Mm$ is representable. Since $\RA$ is a discriminator variety,
$\Mm$ is a subdirect product of some simple algebras $\Cc_i$. We
will show that each $\Cc_i\in\RRA$, this will imply that
$\Mm\in\RRA$.
Now, $\Cc_i$ is in $\Var(\RRA^c)$ since it is a homomorphic image of
$\Mm\in\Var(\RRA^c)$.  Each $\Cc_i$ remains measurable, this is easy
to check directly by using the definition of a measurable relation
algebra, but this fact also follows from \cite[Theorem 3.1]{gajsl}
and \cite[Theorems 6.1, 6.2]{andgiv1}.
Thus, each $\Cc_i$ is persistently finite in $\RA$ by Theorem
\ref{finmeas-thm}. Thus, we can use Theorem \ref{discr-thm} with
substituting $\RA$, $\RRA$ and $\Cc_i$ in place of $\V$, $\K$ and
$\Aa$ and using that $\RA$ is an ordered discriminator variety and
$\Var(\RRA^c)=\Var(\RRA^d)$. We get that $\Cc_i\in\Var\RRA=\RRA$, as
we wanted, and so $\Mm\in\RRA$ since it is a subdirect product of
the $\Cc_i$s.

(ii) follows from (iv) by using the main result of \cite{AGN18} that
there are infinitely many simple finite non-representable measurable
relation algebras. All these algebras are non-integral. This can be
seen by looking at the construction, or by noticing that simple
integral measurable relation algebras are functionally dense and
hence representable.

The proof of (i) is similar to that of (ii): in place of the
non-representable persistently finite measurable relation algebras
of \cite{AGN18} we use the infinitely many non-representable
integral simple maximal relation algebras constructed in \cite{FM}.

To prove (iii), we can use any infinite sequence $\langle \Aa_i :
i\in I\rangle$ of simple, persistently finite non-representable
relation algebras. We have seen in the proofs of (ii) and (i) that
there is such a sequence. We have seen that each persistently finite
relation algebra has only finitely many simple extensions (see the
proof of (5) in the proof of Theorem \ref{discr-thm}). Therefore,
there is an infinite sub-sequence $\langle \Aa_j : j\in J\rangle$ of
the original one such that no $\Aa_j$ can be embedded into $\Aa_k$
for distinct $j,k\in J$. Having this sequence, we will repeat the
proof of \cite[Theorem 5.1]{AGN18} with the necessary modifications.
For $S\subseteq J$ let
\[ \V(S) = \{ \Bb\in\RA : \Aa_n\mbox{ cannot be embedded into }\Bb,
\mbox{ for all }n\in S\}.\] First we show that $\V(S)$ is a variety.
By \cite[Theorem 7.1]{j82}, each finite simple relation algebra is
splitting in the class $\RA$ of all relation algebras, thus the
biggest variety of $\RA$ not containing $\Aa_n$ is the class of all
relation algebras into which $\Aa_n$ cannot be embedded.  This is
called the conjugate variety of $\Aa_n$, let us denote it by
$\V^-(\Aa_n)$, then  \[ \V(S)=\bigcap\{\V^-(\Aa_n) : n\in S\}.\]
This shows that $\V(S)$ is a variety since it is an intersection of
varieties.

All the $\Aa_n$ are simple, persistently finite and
non-representable, so none of them can be embedded into an element
of $\Var(\RRA^c)$, by Theorem \ref{discr-thm}. Thus
\[ \Var(\RRA^c)\subseteq\V(S)\subseteq\RA.\]
We are going to show that $\V(S)$ is distinct from $\V(Z)$ for
distinct subsets $S,Z$ of $J$. This will suffice, because $J$ is
countably infinite, and so it has continuum many subsets. Indeed,
let $S,Z$ be distinct subsets of $J$. Then there is an $n\in J$ such
that, say, $n\in S$ and $n\notin Z$. Then $\Aa_n\notin\V(S)$ since
it can be embedded into itself and $n\in S$. On the other hand,
$\Aa_n\in\V(Z)$ because $n\notin Z$, so no $\Aa_m$ with $m\in Z$ can
be embedded into $\Aa_n$. This shows that $\V(S)\neq\V(Z)$ and we
are done with proving (iii). \qed

The splittable relation algebras, see \cite[Definition 4]{amadn},
all fail to be persistently finite, because once an algebra has a
splittable atom, this atom can be split to arbitrarily many parts
(see \cite[Theorem 3]{amadn}). On the other hand, the
non-representable relation algebras that are shown to be in
$\Var(\RRA^c)$ are gotten by splitting finite simple relation
algebras \cite{maddux}. Perhaps, the next question to ask is the
following.

\begin{prb}\label{split-in} Is every finite simple splittable
relation algebra in $\Var(\RRA^c)$?
\end{prb}

\bigskip\bigskip\bigskip

\noindent Alfr\'ed R\'enyi Institute of Mathematics, Hungarian
Academy of Sciences\\
Budapest, Re\'altanoda st.\ 13-15, H-1053 Hungary\\
andreka.hajnal@renyi.mta.hu, nemeti.istvan@renyi.mta.hu

\end{document}